\documentclass{ifacconf}

\usepackage{graphicx}      % include this line if your document contains figures
\usepackage{amsmath} 	% assumes amsmath package installed
\usepackage{amssymb}  	% assumes amsmath package installed
\usepackage{amsfonts}

\usepackage{subfigure}
\usepackage{epstopdf}
\usepackage{dblfloatfix}
\usepackage{array,xcolor}
\usepackage{multirow}
\usepackage{booktabs}

\usepackage{natbib}        % required for bibliography

\usepackage{url}
\usepackage{balance}

\newtheorem{theorem}{Theorem}

\begin{document}
\begin{frontmatter}

\title{Fast ADMM for homogeneous self-dual embedding of sparse SDPs\thanksref{footnoteinfo}}

\thanks[footnoteinfo]{Y.~Zheng and G.~Fantuzzi contributed equally to this work. Y. Zheng is supported by the Clarendon Scholarship and the Jason Hu Scholarship.}

\author[First]{Yang~Zheng}
\author[Second]{Giovanni~Fantuzzi}
\author[First]{Antonis~Papachristodoulou}
\author[First]{Paul~Goulart}
\author[Second]{and~Andrew~Wynn}

\address[First]{Department of Engineering Science, University of Oxford, Parks Road, Oxford, OX1 3PJ, U.K. (e-mail:
yang.zheng@eng.ox.ac.uk; paul.goulart@eng.ox.ac.uk; antonis@eng.ox.ac.uk).}
\address[Second]{Department of Aeronautics, Imperial College London, South Kensington Campus, London, SW7 2AZ, U.K. (e-mail: gf910@ic.ac.uk; a.wynn@imperial.ac.uk).}

\begin{abstract}
% Abstract of not more than 250 words.
We propose an efficient first-order method, based on the alternating direction method of multipliers (ADMM), to solve the homogeneous self-dual embedding problem for a primal-dual pair of semidefinite programs (SDPs) with chordal sparsity. Using a series of block eliminations, the per-iteration cost of our method is the same as applying a splitting method to the primal or dual alone. Moreover, our approach is more efficient than other first-order methods for generic sparse conic programs since we work with smaller semidefinite cones. In contrast to previous first-order methods that exploit chordal sparsity, our algorithm returns both primal and dual solutions when available, and a certificate of infeasibility otherwise. Our techniques are implemented in the open-source MATLAB solver CDCS. Numerical experiments on three sets of benchmark problems from the library SDPLIB show speed-ups compared to some common state-of-the-art software packages.
\end{abstract}


\end{frontmatter}
%===============================================================================

\section{Introduction}

Semidefinite programs (SDPs) are convex optimization problems commonly used in control theory, machine learning and signal processing. It is well known that although small and medium-sized SDPs can be efficiently solved in polynomial time using second-order interior-point methods (IPMs), these methods become less practical for large-scale SDPs due to memory and time constraints~\citep{helmberg1996interior}. As noted by~\cite{andersen2011interior}, exploiting sparsity in SDPs has been one of the main approaches to improve the scalability of semidefinite programming, and it is still an active and challenging area of research.

In this paper, we present an efficient first-order algorithm to solve the homogeneous self-dual embedding formulation of large-scale SDPs characterized by \textit{chordal sparsity}, meaning that the graph representing their aggregate sparsity pattern is chordal (or has a sparse \textit{chordal extension}). Chordal graphs---undirected graphs with the property that every cycle of length greater than three has a chord---are very well studied objects in graph theory~\citep{vandenberghe2014chordal}. Their connection to SDPs relies on two fundamental theorems due to~\cite{grone1984positive} and~\cite{agler1988positive}: provided that its sparsity pattern is chordal, a large positive semidefinite (PSD) cone can be equivalently replaced with a set of smaller PSD cones.

For this reason chordal sparsity is a key feature of SDPs, %~\citep{de2010exploiting},
and recent years have seen increasing efforts to exploit it in order to increase the computational efficiency of SDP solvers.
%~\citep{vandenberghe2014chordal, zheng2016fast}.
For instance, \cite{fukuda2001exploiting} and \cite{kim2011exploiting} proposed the \emph{domain-space} and the \emph{range-space} conversion techniques to reduce the computational burden of existing IPMs for sparse SDPs. These techniques, implemented in the MATLAB package SparseCoLO~\citep{fujisawa2009user}, rely on the introduction of additional equality constraints to decouple the smaller PSD cones obtained from Grone's and Agler's theorems. However, the addition of equality constraints often offsets the benefit of working with smaller semidefinite cones.

One possible solution to this problem is to exploit the properties of chordal sparsity directly in IPMs~\citep{fukuda2001exploiting, andersen2010implementation}. Another promising direction is to solve decomposable SDPs via first-order methods. For instance, \cite{sun2014decomposition} proposed a first-order splitting algorithm for conic optimization with partially separable structure. \cite{Kalbat2015Fast} applied the alternating direction method of multipliers (ADMM) to solve a special class of SDPs with fully decomposable constraints. \cite{Madani2015ADMM} developed a highly-parallelizable ADMM algorithm for sparse SDPs with applications to optimal power flow problems. More recently, the authors have combined ADMM and chordal decomposition to solve sparse SDPs in either primal or dual standard forms~\citep{zheng2016fast}, providing a conversion framework which is  suitable for the application of first-order methods and parallels that of~\cite{fukuda2001exploiting} and~\cite{kim2011exploiting} for IPMs.

However, none of the aforementioned first-order methods can handle infeasible or unbounded problems. Solving the homogeneous self-dual embedding of the primal-dual pair of optimization problems~\citep{ye1994nl} provides an elegant solution to this issue. The essence of this method is to search for a non-zero point in the non-empty intersection of a convex cone and an affine space. Using this point, one can then either recover an optimal solution of the original primal-dual pair of SDPs, or construct a certificate of primal or dual infeasibility. Homogeneous self-dual embeddings have been widely used in IPMs~\citep{sturm1999using,ye2011interior}; more recently, \cite{ODonoghue2016} have proposed an operator-splitting method for the homogeneous self-dual embedding of general conic programs that scales well with problem size, and the implementation can be found in the C package \textsc{SCS}~\citep{scs}.

In this work, we show that the conversion techniques for primal and dual standard-form SDPs developed in~\cite{zheng2016fast} can be extended to the homogeneous self-dual embedding. Also, we extend the algorithm in~\cite{ODonoghue2016} to take advantage of chordal sparsity. Our main contributions are:

\begin{enumerate}
  \item We formulate the homogeneous self-dual embedding of a primal-dual pair of SDPs whose conic constraints are decomposed using Grone's and Agler's theorems. This extends the conversion techniques for sparse SDPs developed in our previous work~\citep{zheng2016fast}.  To the best of our knowledge, it is the first time that such a formulation has been presented.
  \item  We extend the ADMM algorithm of \cite{ODonoghue2016} to take advantage of the special structure of our homogeneous self-dual formulation, thereby reducing its computational complexity. Our algorithm is more efficient than a direct application of the method of \cite{ODonoghue2016} to either the original primal-dual pair (i.e. before chordal sparsity is taken into account), or the decomposed problems: in the former case, the chordal decomposition reduces the cost of the conic projections; in the latter case, we speed up the affine projection step using a series of block-eliminations.
  \item We implement our techniques in the MATLAB solver CDCS (Cone Decomposition Conic Solver). This is the first open source first-order solver that exploits chordal decomposition and is able to handle infeasible problems. Numerical simulations on three sets of benchmark problems from the library SDPLIB~\citep{borchers1999sdplib} demonstrate the efficiency of our self-dual algorithm compared to other commonly used software packages.
\end{enumerate}

The rest of this paper is organized as follows. Section~\ref{se:preliminaries} reviews some background material. We present the homogeneous self-dual embedding of SDPs with chordal sparsity in Section~\ref{se:HomogeneousEmbedding}. Section~\ref{se:ADMM} discusses our ADMM algorithm in detail, and we report numerical experiments in Section~\ref{se:simulation}. Finally, Section~\ref{se:conclusion} offers concluding remarks.

\section{Preliminaries} \label{se:preliminaries}

\subsection{Chordal graphs}
Let $\mathcal{G}(\mathcal{V},\mathcal{E})$ be an undirected graph with nodes $\mathcal{V}=\{1,2,\ldots,n\}$ and edges $\mathcal{E} \subseteq \mathcal{V} \times \mathcal{V}$. A subset of nodes $\mathcal{C}\subseteq \mathcal{V} $ is called a \emph{clique} if $ (i,j) \in \mathcal{E}$ for any distinct nodes $ i,j \in \mathcal{C} $.  If $\mathcal{C}$ is not a subset of any other clique, then it is called a \emph{maximal clique}.  The number of nodes in $\mathcal{C}$ is denoted by $\vert \mathcal{C} \vert$, and $\mathcal{C}(i)$ indicates the $i$-th element of $\mathcal{C}$, sorted in the natural ordering.
An undirected graph $\mathcal{G}$ is called \emph{chordal} if every cycle of length greater than three has at least one chord (an edge connecting two nonconsecutive nodes in the cycle). Note that if $\mathcal{G}(\mathcal{V}, \mathcal{E})$ is not chordal, it can be \emph{chordal extended}, \emph{i.e.}, we can construct a chordal graph $\mathcal{G}'(\mathcal{V}, \mathcal{E}') $ by adding suitable edges to $ \mathcal{E} $.%~\citep{yannakakis1981computing}.

\subsection{Sparse matrices defined by graphs}

Let $ \mathcal{G}(\mathcal{V},\mathcal{E}) $ be an undirected graph, and assume that $(i,i) \in \mathcal{E}$ for any node $i \in \mathcal{V}$. A partial symmetric matrix is a symmetric matrix in which the entry $X_{ij}$ is specified if and only if $(i,j) \in \mathcal{E}$. In this work, we use the following sets of symmetric matrices defined on $\mathcal{E}$:
\begin{align*}
    \mathbb{S}^n(\mathcal{E},?) = &\text{ the space of } n \times n \text{ partial symmetric matrices} \\[-.25em]
                                  &\text{ with elements defined on } \mathcal{E}, \\
    \mathbb{S}_{+}^n(\mathcal{E},?) = &\{ X \in \mathbb{S}^n(\mathcal{E},?) \mid \exists M \succeq 0, M_{ij} = X_{ij}, \forall (i,j) \in {\mathcal{E}}  \}, \\
    \mathbb{S}^n(\mathcal{E},0) = &\{ X \in \mathbb{S}^n \mid X_{ij} = 0, \text{if } (i,j) \notin {\mathcal{E}}  \}, \\
    \mathbb{S}_{+}^n(\mathcal{E},0) = &\{ X \in \mathbb{S}^n(\mathcal{E},0) \mid X \succeq 0 \}.
\end{align*}

Note that $ \mathbb{S}^n_{+}(\mathcal{E},?)$ and $\mathbb{S}_{+}^n(\mathcal{E},0)$ are two types of sparse matrix cones, and that they are the dual of each other for any (that is, chordal or not) sparsity pattern $\mathcal{E}$~\citep{vandenberghe2014chordal}.

Finally, let $\mathcal{C}$ be a maximal clique of the graph $\mathcal{G}$, and let $E_{\mathcal{C}} \in \mathbb{R}^{\mid \mathcal{C}\mid \times n}$ be the matrix with entries $(E_{\mathcal{C}})_{ij} = 1$ if $\mathcal{C}(i) = j$ and $(E_{\mathcal{C}})_{ij} = 0$ otherwise. Then, given a symmetric matrix $X \in \mathbb{S}^n$, the submatrix of $X$ defined by the clique $\mathcal{C}$ can be represented as $E_{\mathcal{C}}XE_{\mathcal{C}}^T \in \mathbb{S}^{\mid \mathcal{C}\mid}$.

\subsection{Chordal decomposition}

The problems of deciding if $X \in \mathbb{S}^n_{+}(\mathcal{E},?)$ or $Z \in \mathbb{S}_{+}^n(\mathcal{E},0)$ can be posed as problems over several smaller (but coupled) convex cones according to the following theorems:
\begin{theorem}\label{T:ChordalCompletionTheorem}\citep{grone1984positive}.
     Let $\mathcal{G}(\mathcal{V},\mathcal{E})$ be a chordal graph with maximal cliques $\{\mathcal{C}_1,\mathcal{C}_2, \ldots, \mathcal{C}_p\}$. Then, $X\in\mathbb{S}^n_+(\mathcal{E},?)$ if and only if $X_k := E_{\mathcal{C}_k} X E_{\mathcal{C}_k}^T \in \mathbb{S}^{\vert \mathcal{C}_k \vert}_+$ for all $k=1,\,\ldots,\,p$.
\end{theorem}

\begin{theorem}\label{T:ChordalDecompositionTheorem}\citep{agler1988positive}.
     Let $\mathcal{G}(\mathcal{V},\mathcal{E})$ be a chordal graph with maximal cliques $\{\mathcal{C}_1,\mathcal{C}_2, \ldots, \mathcal{C}_p\}$. Then, $Z\in\mathbb{S}^n_+(\mathcal{E},0)$ if and only if there exist matrices $Z_k \in \mathbb{S}^{\vert \mathcal{C}_k \vert}_+$ for $k=1,\,\ldots,\,p$ such that $Z = \sum_{k=1}^{p} E_{\mathcal{C}_k}^T Z_k E_{\mathcal{C}_k}.$
\end{theorem}

%Note that these two theorems can be proven individually, but can also can be derived from each other using the duality of the cones $\mathbb{S}^n_{+}(\mathcal{E},?)$ and $\mathbb{S}_{+}^n(\mathcal{E},0)$~\citep{vandenberghe2014chordal}.

\section{Homogeneous Self-Dual Embedding of Sparse SDPs} \label{se:HomogeneousEmbedding}

% Is the following really needed? Maybe better to just go straight to the point
% In this section we present the homogeneous self-dual embedding of sparse SDPs, in which Grone's and Agler's theorems are used to decompose the primal and dual positive semidefinite constraints. The presentation in Section~\ref{SS:SparseChordalDec} closely follows that of~\cite{zheng2016fast}, and has been included to make the current work self-contained.

Consider the standard \emph{primal-dual pair} of SDPs, \emph{i.e.},
\begin{equation} \label{E:PrimalSDP}
    \begin{aligned}
    \min_{X} \quad & \langle C, X \rangle \\
    \text{subject to} \quad & \langle A_i,X \rangle = b_i,
    \quad i =1 ,\ldots, m,\\
        & X \in \mathbb{S}^n_{+},
    \end{aligned}
\end{equation}
and
\begin{equation} \label{E:DualSDP}
    \begin{aligned}
        \max_{y,Z} \quad & \langle b,y \rangle\\
        \text{subject to} \quad  & \sum_{y=1}^m A_i y_i + Z= C,\\
        & Z \in \mathbb{S}^n_{+}.
    \end{aligned}
\end{equation}
The vector $ b\in \mathbb{R}^m$ and  the matrices $C,\,A_1,\,\ldots,\,A_m \in \mathbb{S}^n$ are the problem data; $X$ is the primal variable, and $y,Z$ are the dual variables. We say that \eqref{E:PrimalSDP} and \eqref{E:DualSDP} have the \emph{aggregate sparsity pattern} $\mathcal{G}(\mathcal{V},\mathcal{E})$ if
$C,A_1,\ldots,A_m \in \mathbb{S}^{n}(\mathcal{E},0)$.
Throughout this work, we will assume that $\mathcal{G}$ is chordal (otherwise, it can be chordal extended), and that its maximal cliques $\mathcal{C}_1,\ldots,\mathcal{C}_p$ are small.

\subsection{Sparse SDPs with Chordal Decomposition}
\label{SS:SparseChordalDec}
Aggregate sparsity implies that the dual variable $Z$ in \eqref{E:DualSDP} must have the sparsity pattern defined by $\mathcal{E}$, \emph{i.e.}, $Z \in \mathbb{S}^{n}(\mathcal{E},0)$. Similarly, although the primal variable $X$ in \eqref{E:PrimalSDP} is usually dense, the cost function and the equality constraints only depend on the entries $X_{ij}$ in the sparsity pattern $\mathcal{E}$, while the remaining entries only guarantee that $X$ is positive semidefinite. This means that it suffices to consider $X\in \mathbb{S}^n_{+}(\mathcal{E},?)$. Then, according to Theorems~\ref{T:ChordalCompletionTheorem}-\ref{T:ChordalDecompositionTheorem}, we can rewrite \eqref{E:PrimalSDP} and \eqref{E:DualSDP}, respectively, as
\begin{equation}\label{E:DecomposedPriaml}
  \begin{aligned}
    \min_{X,X_1,\ldots,X_p} \quad & \langle C, X \rangle \\
    \text{subject to} \quad  & \langle A_i, X \rangle = b_i, &i &= 1, \ldots, m,\\
        & X_k = E_{\mathcal{C}_k} X E_{\mathcal{C}_k}^T, &k &= 1, \ldots, p, \\
        & X_k \in \mathbb{S}^{\mid \mathcal{C}_k\mid}_{+}, &k &= 1, \ldots, p,
    \end{aligned}
\end{equation}
and
\begin{equation}\label{E:DecomposedDual}
  \begin{aligned}
        \max_{y,Z_1,\ldots,Z_p,V_1,\ldots,V_p} \quad & \langle b,y \rangle \\
        \text{subject to} \quad &  \sum_{i=1}^m A_i y_i + \sum_{k=1}^p E_{\mathcal{C}_k}^T V_k E_{\mathcal{C}_k} = C,\\
        & Z_k = V_k, \qquad k = 1, \ldots, p,\\
        & Z_k \in \mathbb{S}^{\mid \mathcal{C}_k\mid}_{+}, \quad\,\, k = 1, \ldots, p.
    \end{aligned}
\end{equation}

%It is not difficult to check that the decomposed problems \eqref{E:DecomposedPriaml} and \eqref{E:DecomposedDual} are also the dual of each other by virtue of the duality between Grone's and Agler's theorems.

To ease the exposition, let $\mathrm{vec}:\mathbb{S}^n \to \mathbb{R}^{n^2}$ be the usual operator mapping a matrix to the stack of its column, and define the vectorized data $ c := \text{vec}(C)$,  $A := [\text{vec}(A_0)\,\hdots\, \text{vec}(A_m)]^T$,
%
%\begin{align*}
%    c &:= \text{vec}(C), & A := \begin{bmatrix} \text{vec}(A_0) & \hdots & \text{vec}(A_m) \end{bmatrix}^T,
%\end{align*}
the vectorized variables
%\begin{align*}
        $x := \mathrm{vec}(X)$,
        $x_k := \text{vec}(X_k)$,
        $z_k := \text{vec}(Z_k)$,
        $v_k := \text{vec}(V_k)$   for all $k = 1,\,\ldots,\,p$,
%\end{align*}
%
%\begin{align*}
%        x &:= \mathrm{vec}(X),  &
%        x_k &:= \text{vec}(X_k),   &  \\
%        z_k &:= \text{vec}(Z_k),  &
%        v_k &:= \text{vec}(V_k),   & k &= 1,\,\ldots,\,p,
%\end{align*}
and the matrices
\begin{equation}
    H_k := E_{\mathcal{C}_k} \otimes E_{\mathcal{C}_k},\quad k = 1, \ldots, p,
\end{equation}
such that
%\begin{equation*}
    $x_k = \mathrm{vec}(X_k) = \mathrm{vec}(E_k X E_k^T) = H_k x$.
%\end{equation*}
Note that $H_1,\,\ldots,\,H_p$ are ``entry-selector'' matrices of $1$'s  and $0$'s, whose rows are orthonormal and such that $H_k^TH_k$ is diagonal. These matrices project $x$ onto the subvectors $x_1,\,\ldots,\,x_p$, respectively.

If we denote the constraints $X_k\in\mathbb{S}^{|\mathcal{C}_k|}_+$ by $x_k\in\mathcal{S}_k$, we can rewrite~\eqref{E:DecomposedPriaml} and~\eqref{E:DecomposedDual} as
\begin{equation} \label{E:PrimalVectorForm}
\begin{aligned}
\min_{x,x_1,\ldots, x_p} \quad & \langle c,x \rangle\\
\text{subject to} \quad & Ax=b,\\
					   & x_k = H_k x,&& k=1,\,\ldots,\,p,\\
					   & x_k \in \mathcal{S}_k,&& k=1,\,\ldots,\,p,
\end{aligned}
\end{equation}
and
\begin{equation}\label{E:DualVectorForm}
    \begin{aligned}
        \max_{y,z_1,\ldots,z_p,v_1,\ldots,v_p} \quad & \langle b,y \rangle\\
        \text{subject to} \quad
        & A^Ty + \sum_{k=1}^p H_k^Tv_k =  c ,\\
        & z_k - v_k = 0, \quad k=1,\,\ldots,\,p,\\
        & z_k \in \mathcal{S}_k, \qquad\,\,\,\, k=1,\,\ldots,\,p.
    \end{aligned}
\end{equation}

\subsection{Homogeneous Self-Dual Embedding}

For notational simplicity, let $\mathcal{S} := \mathcal{S}_1 \times \cdots \times \mathcal{S}_p$ and define
\begin{align*}
  s &:= \begin{bmatrix}x_1 \\ \vdots\\ x_p \end{bmatrix}, &
        z &:= \begin{bmatrix}z_1 \\ \vdots\\ z_p \end{bmatrix}, &
        v &:= \begin{bmatrix}v_1 \\ \vdots\\ v_p \end{bmatrix}, &
        H &:= \begin{bmatrix}H_1 \\ \vdots\\ H_p \end{bmatrix}.
\end{align*}

When strong duality holds for \eqref{E:PrimalVectorForm} and \eqref{E:DualVectorForm}, the following KKT conditions are necessary and sufficient for optimality of the tuple
%$(x^*,r^*,w^*, s^*,h^*,y^*,v^*,z^*)$:
$(x,s,y,v,z)$:
\begin{itemize}
  \item  %$(x^*,r^*,w^*, s^*)$
  		$(x, s)$ is primal feasible; adding optimal slack variables
  		$r=0$ and $w=0$, this amounts to
        \begin{equation}\label{E:PrimalFeasible}
            \begin{aligned}
                 Ax - r &= b,  & r &= 0, \\
                s + w &= Hx,  & w &= 0, & s &\in \mathcal{S}.
            \end{aligned}
        \end{equation}
  \item %$(h^*,y^*,v^*,z^*)$
  		$(y,v,z)$	is dual feasible; adding an optimal slack variable
  		$h=0$, this amounts to
        \begin{equation}\label{E:DualFeasible}
            \begin{aligned}
                 A^Ty + H^Tv + h &=  c , & h &= 0, \\
                z - v &= 0, & z &\in \mathcal{S}.
            \end{aligned}
        \end{equation}
  \item The duality gap is zero,  \textit{i.e.}
        \begin{equation}\label{E:ZeroGap}
            c^Tx - b^Ty = 0.
        \end{equation}
\end{itemize}

The idea behind the homogeneous self-dual embedding~\citep{ye1994nl} is to introduce two non-negative and complementary variables $\tau$ and $\kappa$ in addition to the variables zero variables $r$, $w$ and $h$ introduced above, and embed the KKT conditions \eqref{E:PrimalFeasible}, \eqref{E:DualFeasible} and \eqref{E:ZeroGap} into the linear system
\begin{equation} \label{E:HomogeneousEmbedding}
 \begin{bmatrix}
    h \\ z \\ r \\ w \\ \kappa
  \end{bmatrix} = \begin{bmatrix}
    0 & 0 & -A^T & -H^T & c \\
    0 & 0 & 0 & I & 0\\
    A & 0 & 0 & 0 & -b\\
    H & -I & 0 & 0 & 0 \\
    -c^T & 0  & b^T & 0 & 0\\
  \end{bmatrix} \begin{bmatrix}
    x \\ s \\ y \\ v \\ \tau
  \end{bmatrix}.
\end{equation}
Any solution of this embedding can be used to recover an optimal solution for \eqref{E:PrimalVectorForm}-\eqref{E:DualVectorForm}, or provide a certificate of primal or dual infeasibility~\citep[see][]{ODonoghue2016}.

Letting $n_d =\sum_{k=1}^p |\mathcal{C}_k|^2$, defining $\mathcal{K} := \mathbb{R}^{n^2} \times \mathcal{S} \times \mathbb{R}^{m} \times \mathbb{R}^{n_d} \times \mathbb{R}_{+}$,
%
%\[
%\mathcal{K} := \mathbb{R}^{n^2} \times \mathcal{S} \times \mathbb{R}^{m} \times \mathbb{R}^{n_d} \times \mathbb{R}_{+},
%\]
and writing
\begin{align*}
  u &:= \begin{bmatrix}
    x \\ s \\ y \\ v \\ \tau
  \end{bmatrix},
  &
  v &:= \begin{bmatrix}
    h \\ z \\ r \\ w \\ \kappa
  \end{bmatrix},
  &
  Q &:= \begin{bmatrix}
    0 & 0 & -A^T & -H^T & c \\
    0 & 0 & 0 & I & 0\\
    A & 0 & 0 & 0 & -b\\
    H & -I & 0 & 0 & 0 \\
    -c^T & 0  & b^T & 0 & 0\\
  \end{bmatrix}
\end{align*}
to further ease the notation, the decomposed primal-dual pair of SDPs~\eqref{E:PrimalVectorForm}-\eqref{E:DualVectorForm} can be recast as the feasibility problem
\begin{equation} \label{E:ADMMForm}
\begin{aligned}
\text{find} \qquad &(u,v)\\
\text{subject to} \quad &v = Qu,\\
			      &(u,v) \in \mathcal{K} \times \mathcal{K}^*,
\end{aligned}
\end{equation}

where $\mathcal{K}^*=\{0\}^{n^2}\!\times\mathcal{S}^*\!\times\!\{0\}^{m}\!\times\!\{0\}^{n_d}\!\times\!\mathbb{R}_{+}$ is the dual of the cone $\mathcal{K}$ (here, $\mathcal{S}^*$ is the dual of $\mathcal{S}$, and $\{0\}^p$ is the $p$-dimensional zero vector).

\section{ADMM for the Homogeneous Self-dual Embedding} \label{se:ADMM}
\label{S:ADMM}

\subsection{Basic algorithm}

Problem~\eqref{E:ADMMForm} is in the same form considered by \cite{ODonoghue2016}, so we can directly apply their ADMM algorithm. The $k$-th iteration of the algorithm consists of the following three steps \citep{ODonoghue2016}, where $\Pi_{\mathcal{K}}$ denotes projection on the cone $\mathcal K$:
\begin{subequations} \label{E:ADMMSteps}
    \begin{align}
    \label{E:Step1}
      \hat{u}^{k+1} &= (I+Q)^{-1}(u^k+ v^k), \\
    \label{E:Step2}
      u^{k+1} &= \Pi_{\mathcal{K}}(\hat{u}^{k+1}-v^k),\\
    \label{E:Step3}
      v^{k+1} &= v^k - \hat{u}^{k+1} + u^{k+1}.
    \end{align}
\end{subequations}

Note that~\eqref{E:Step2} is inexpensive, since $\mathcal{K}$ is the cartesian product of simple cones (zero, free and non-negative cones) and small PSD cones, and can be efficiently carried out in parallel. The third step is also computationally inexpensive and parallelizable. On the contrary, although the preferred factorization of $I+Q$ (or its inverse) can be cached before starting the iterations, a direct implementation of~\eqref{E:Step1} can be computationally intensive since $Q$ is a very large matrix. However, $Q$ is highly structured and sparse; in the next sections, we show how its special structure can be exploited to devise a series of block-eliminations that speed up the affine projection in~\eqref{E:Step1}.

\subsection{Solving the linear system}

The affine projection step~\eqref{E:Step1} requires the solution of a linear system in the form
\begin{equation} \label{E:LinearSystem}
\begin{bmatrix}
I & \hat{A}^T & \hat{c} \\
-\hat{A} & I & \hat{b} \\
-\hat{c}^T & -\hat{b}^T & 1
\end{bmatrix}
\begin{bmatrix}
\hat{u}_1 \\ \hat{u}_2 \\ \hat{u}_3
\end{bmatrix}
=
\begin{bmatrix}
\omega_1 \\ \omega_2 \\ \omega_3
\end{bmatrix},
\end{equation}
where
\begin{align*}
\hat{A} &= \begin{bmatrix} -A & 0 \\ -H & I \end{bmatrix}, &
\hat{c} &= \begin{bmatrix}c\\0 \end{bmatrix},&
\hat{b} &= \begin{bmatrix}-b\\0 \end{bmatrix}.
\end{align*}
Note that $u_3$ and $\omega_3$ are scalars. Letting
\begin{align*}
M &:= \begin{bmatrix}
I & \hat{A}^T \\
-\hat{A} & I
\end{bmatrix},
&
\zeta &:= \begin{bmatrix}\hat{c} \\ \hat{b} \end{bmatrix},
\end{align*}
and carrying out block elimination on \eqref{E:LinearSystem} we obtain
\begin{subequations}
\begin{align}
(M+\zeta\zeta^T) \begin{bmatrix}\hat{u}_1 \\ \hat{u}_2  \end{bmatrix} &= \begin{bmatrix}\omega_1 \\ \omega_2  \end{bmatrix} - \omega_3 \zeta. \label{E:MatrixInverse} \\
\hat{u}_3 &= \omega_3 + \hat{c}^T \hat{u}_1 + \hat{b}^T\hat{u}_2. \label{E:FirstBlockElimination}
\end{align}
\end{subequations}
Moreover, the matrix inversion lemma~\citep[Appendix C.4.3]{boyd2004convex} allows us to write the solution of \eqref{E:MatrixInverse} as
\begin{equation} \label{E:MatrixInverseResult}
    \begin{bmatrix}\hat{u}_1 \\ \hat{u}_2  \end{bmatrix}  =
    \left[
    M^{-1} -
    \frac{ (M^{-1} \zeta)\zeta^T M^{-1}}{1 + \zeta^T (M^{-1} \zeta)}
    \right]
    \left(
    \begin{bmatrix}\omega_1 \\ \omega_2  \end{bmatrix} - \omega_3 \zeta
    \right).
    %\\
    %&=
    %M^{-1} \left( v - v_3 \zeta \right) -
    %\frac{ (M^{-1} \zeta) \zeta^T }{1 + \zeta^T (M^{-1} \zeta)}
    %\left[  M^{-1} \left( v - v_3 \zeta \right)\right]
\end{equation}
Note that the vector $M^{-1} \zeta$ and the scalar $1 + \zeta^T (M^{-1} \zeta)$ only depend on the problem data, and can be cached before starting the ADMM iterations. %Instead, $M^{-1} \left( v - v_3 \zeta \right)$ must be computed at each iteration since $v - v_3 \zeta$ changes at each iteration.
Consequently, updating $\hat{u}_1$, $\hat{u}_2$ and $\hat{u}_3$ at each iteration requires:
\begin{enumerate}
\item the solution of the ``inner'' linear system to compute
\[
M^{-1} \left(\begin{bmatrix}\omega_1 \\ \omega_2  \end{bmatrix} - \omega_3 \zeta\right).
\]
\item a series of inexpensive vector inner products and scalar-vector operations in \eqref{E:FirstBlockElimination} and \eqref{E:MatrixInverseResult}.
\end{enumerate}

\begin{table*}[!h]
    \centering
    \renewcommand\arraystretch{0.68}
    \caption{Details of the SDPLIB problems considered in this work.}
    \label{T:SparseStatistic}
    \begin{tabular}{r | c c c c | c c c c | c c}%{m{2cm}<{\centering} m{0.9cm}<{\centering} m{0.9cm}<{\centering} m{0.9cm}<{\centering} m{0.9cm}<{\centering}}
        \hline  \toprule\\[-1.0em]
        & \multicolumn{4}{c|}{Small and medium-size $(n \leq 100)$}  & \multicolumn{4}{c|}{Large-scale and sparse $(n \geq 800)$} & \multicolumn{2}{c}{Infeasible}  \\
            & theta1 & theta2 & qap5 & qap9  & maxG11 & maxG32 &  qpG11 & qpG51 & infp1 & infd1   \\
        \hline\\[-0.5em]
        Original cone size, $n$  & 50  & 100 & 26  & 82  & 800 & 2000 & 1600 & 2000 & 30 & 30  \\
        Affine constraints, $m$  & 104 & 498 & 136 & 748 & 800 & 2000 & 800  & 1000 & 10 & 10 \\
        Number of cliques, $p$   & 1   & 1   & 1   & 1   & 598 & 1499 & 1405 & 1675 & 1  & 1  \\
        Maximum clique size      & 50  & 100 & 26  & 82  & 24  & 60   & 24   & 304  & 30 & 30 \\
        Minimum clique size      & 50  & 100 & 26  & 82  & 5   & 5    & 1    & 1    & 30 & 30  \\
        \bottomrule
        \end{tabular}
\end{table*}

\begin{table*}
\centering
\renewcommand\arraystretch{0.68}
\caption{Results for some small and medium-sized SDPs in SDPLIB}
\label{T:ResultsSmall}
\begin{tabular}{l r c c c c c c c}
\toprule[1pt]
& & SeDuMi &\begin{tabular}[x]{@{}c@{}}SparseCoLO+\\SeDuMi\end{tabular}
& SCS &\begin{tabular}[x]{@{}c@{}}CDCS\\(primal)\end{tabular}
&\begin{tabular}[x]{@{}c@{}}CDCS\\(dual)\end{tabular}
&\begin{tabular}[x]{@{}c@{}}Self-dual \end{tabular}\\
\midrule
\multirow{4}{*}{theta1} &
Total time (s) 	  & 0.262 & 0.279 & 0.145 & 0.751  & 0.707 & 0.534\\
& Pre- time (s) &  0    & 0.005 & 0.011  & 0.013 & 0.010 & 0.012\\
& Iterations    & 14 & 14 & 240 & 317 & 320  & 230 \\
& Objective     & 2.300 $\times 10^1$ & 2.300$\times 10^1$ & 2.300$\times 10^1$ & 2.299$\times 10^1$  & 2.299$\times 10^1$ & 2.303$\times 10^1$  \\
\midrule
\multirow{4}{*}{theta2}&
Total time (s)   & 1.45 & 1.55 & 0.92 & 1.45    & 1.30  & 0.60 \\
& Pre- time (s)&  0   & 0.014 & 0.018  & 0.046 & 0.036 & 0.031\\
& Iterations      & 15 & 15 & 500 & 287   & 277  & 110 \\
& Objective       & 3.288 $\times 10^1$ & 3.288$\times 10^1$ & 3.288$\times 10^1$ & 3.288$\times 10^1$ & 3.288$\times 10^1$ & 3.287$\times 10^1$ \\
\midrule
\multirow{4}{*}{qap5}  &
Total time (s)    & 0.365 & 0.386 & 0.412 &  0.879 & 0.748 & 1.465\\
& Pre- time (s)   &  0   & 0.006 & 0.026  &  0.011 & 0.009 & 0.009 \\
& Iterations       & 12 & 12 & 320 & 334 & 332  & 783 \\
& Objective        & -4.360$\times 10^2$ & -4.360$\times 10^2$ & -4.359$\times 10^2$ & -4.360$\times 10^2$  & -4.364$\times 10^2$  & -4.362$\times 10^2$ \\
\midrule
\multirow{4}{*}{qap9}  &
Total time (s)           & 6.291 & 6.751 & 3.261 &  7.520 & 7.397 & 1.173\\
& Pre- time (s)&  0   & 0.012 & 0.010  &  0.064 & 0.036 & 0.032 \\
& Iterations             & 25 & 25 & 2000 & 2000 & 2000  & 261 \\
& Objective        & -1.410$\times 10^3$ & -1.410$\times 10^3$ & -1.409$\times 10^3$ & -1.407$\times 10^3$  & -1.409$\times 10^3$  & -1.410$\times 10^3$ \\
\bottomrule[1pt]
\end{tabular}
\end{table*}

\begin{table*}[!h]
\centering
\renewcommand\arraystretch{0.68}
\caption{Results for some large-scale sparse SDPs}% in SDPLIB.}
\label{T:ResultsLargeScale}
\begin{tabular}{l r c c c c c c c}
\toprule[1pt]
& & SeDuMi &\begin{tabular}[x]{@{}c@{}}SparseCoLO+\\SeDuMi\end{tabular}
& SCS &\begin{tabular}[x]{@{}c@{}}CDCS\\(primal)\end{tabular}
&\begin{tabular}[x]{@{}c@{}}CDCS\\(dual)\end{tabular}
&\begin{tabular}[x]{@{}c@{}}Self-dual \end{tabular}\\
\midrule
\multirow{4}{*}{maxG11} &
Total time (s)   & 92.0 & 9.83 & 160.5 & 126.6   & 114.1  & 23.9 \\
& Pre- time (s)  &  0  & 2.39 & 0.07  & 3.33  & 4.28  & 2.45 \\
& Iterations     & 13 & 15 & 1860 & 1317  & 1306   & 279   \\
& Objective      & 6.292$\times 10^2$ & 6.292$\times 10^2$ & 6.292$\times 10^2$ & 6.292$\times 10^2$  &6.292$\times 10^2$ & 6.295$\times 10^2$  \\
\midrule
\multirow{4}{*}{maxG32} &
Total time (s)   & 1.385 $\times 10^3$ & 577.4 & 2.487 $\times 10^3$ & 520.0   & 273.8 & 87.4 \\
& Pre- time (s)  &  0  & 7.63 & 0.589  & 53.9  & 55.6  & 30.5 \\
& Iterations     & 14 & 15 & 2000 & 1796  & 943   & 272   \\
& Objective      & 1.568$\times 10^3$ & 1.568$\times 10^3$ & 1.568$\times 10^3$ & 1.568$\times 10^3$  &1.568$\times 10^3$ & 1.568$\times 10^3$ \\
\midrule
 \multirow{4}{*}{qpG11} &
Total time (s) 	 & 675.3 & 27.3 & 1.115 $\times 10^3$ & 273.6  & 92.5 & 32.1\\
 & Pre- time (s) & 0    & 11.2 & 0.57  & 6.26 & 6.26 & 3.85\\
 & Iterations    & 14 & 15 & 2000 & 1355 & 656  & 304 \\
 & Objective     & 2.449$\times 10^3$ & 2.449$\times 10^3$ & 2.449$\times 10^3$ & 2.449$\times 10^3$  & 2.449$\times 10^3$ & 2.450$\times 10^3$  \\
\midrule
 \multirow{4}{*}{qpG51}  &
Total time (s)   & 1.984$\times 10^3$ & -- & 2.290$\times 10^3$ &  1.627$\times 10^3$ & 1.635$\times 10^3$ & 538.1\\
& Pre- time (s)  &  0   & -- & 0.90  &  10.82 & 12.77 & 7.89 \\
 & Iterations    & 22 & -- & 2000 & 2000 & 2000  & 716 \\
  & Objective    & 1.182$\times 10^3$ & -- & 1.288$\times 10^3$ & 1.183$\times 10^3$  & 1.186$\times 10^3$  & 1.181$\times 10^3$ \\
\bottomrule[1pt]
\end{tabular}
\end{table*}

\subsection{Solving the ``inner'' linear system}

Recalling the definition of $M$, computing~\eqref{E:MatrixInverseResult} requires the solution of a linear system of the form
\begin{equation} \label{E:SecondBlockElimination}
    \begin{bmatrix}
    I & \hat{A}^T \\
    -\hat{A} & I
    \end{bmatrix}
    \begin{pmatrix}
    \hat{u}_1 \\ \hat{u}_2
    \end{pmatrix}
    = \begin{pmatrix}
    \hat{\omega}_1 \\ \hat{\omega}_2
    \end{pmatrix}.
\end{equation}
Block elimination leads to
\begin{subequations}
\begin{align}
    \hat{u}_1 &= \hat{\omega}_1 - \hat{A}^T \hat{u}_2,\\
    \label{E:SecBlkEliResult}
    (I + \hat{A}\hat{A}^T) \hat{u}_2 &= \hat{A} \hat{\omega}_1 + \hat{\omega}_2.
\end{align}
\end{subequations}
Recalling the definition of $\hat{A}$ and recognizing that
%\[
$D := H^T H= \sum_{k=1}^p H_k^TH_k$
%\]
%\[
%D := H^T H= \sum_{k=1}^p H_k^TH_k
%\]
%
is a diagonal matrix, we also have
\begin{equation} \label{E:ThirdBlkElm}
I + \hat{A}\hat{A}^T = \begin{bmatrix} I+D + A^TA & -H^T \\ -H & 2I\end{bmatrix}.
\end{equation}
Given the special structure of this matrix, block elimination can be used again to solve \eqref{E:SecBlkEliResult}. Simple algebraic manipulations show that the only matrix to be factorized (or inverted) before starting the ADMM iterations is
\begin{equation} \label{E:FactorizationMatrix}
I+\frac{1}{2}\,D + A^T A.
\end{equation}
In fact, this matrix is of the ``diagonal plus low rank" form, so the matrix inversion lemma can be used to reduce the size of the matrix to invert even further.

\subsection{Summary of computational gains}

The algorithm outlined in the previous sections is clearly more efficient than a direct application of the ADMM algorithm of \cite{ODonoghue2016} to the decomposed primal-dual pair of SDPs~\eqref{E:PrimalVectorForm}-\eqref{E:DualVectorForm}. In fact, the cost of the conic projection~\eqref{E:Step2} will be the same for both algorithms, but the sequence of block eliminations and applications of the matrix inversion lemma we have described reduces the cost of the affine projection step. In the algorithm, we only need to factor an $m \times m$ matrix, instead of the matrix $Q$ of size $(n^2 + 2n_d + m+1) \times (n^2 + 2n_d + m+1)$.

Furthermore, it can be checked that when we exploit the special structure of the matrix $I+Q$, the overall computational cost of~\eqref{E:Step1} coincides (at least to leading order) with the cost of the affine projection step when the algorithm of \cite{ODonoghue2016} is applied to the original primal-dual pair~\eqref{E:PrimalSDP}-\eqref{E:DualSDP}, before chordal decomposition.
This means that our algorithm should also outperform the algorithm of \cite{ODonoghue2016} applied to the original primal-dual pair of SDPs~\eqref{E:PrimalSDP}-\eqref{E:DualSDP}: the cost of the affine projection is the same, but the conic projection in our algorithm is more efficient since we work with smaller positive semidefinite cones.

\balance
\section{Numerical Simulations} \label{se:simulation}
We have implemented our techniques in CDCS (Cone Decomposition Conic Solver)~\citep{CDCS}. The code is available at

\begin{center}
\url{https://github.com/oxfordcontrol/CDCS}.
\end{center}
%This is the first open-source first-order solver that exploits chordal decomposition with the ability to handle infeasible problems.
In addition to the homogeneous self-dual embedding algorithm, CDCS also includes the primal and dual methods of~\cite{zheng2016fast}. Currently,  %CDCS supports cartesian products of the following cones: $\mathbb R^n$, non-negative orthant, second-order cone, and positive semidefinite cone. %We only implemented chordal decomposition techniques for semidefinite cones, while the other supported cone types are not decomposed.
CDCS was tested on three sets of benchmark problems in SDPLIB~\citep{borchers1999sdplib}: %\footnote{The experiments reported here are only meant to illustrate the ideas presented in this paper; our codes are under active development, and their behavior might change as improvements are made.}:
{\it(1)} Four small and medium-sized SDPs (two Lov\'asz $\vartheta$ number problems, theta1 and theta2, and two quadratic assignment problems, qap5 and qap9); {\it(2)} Four large-scale sparse SDPs (two max-cut problems, maxG11 and maxG32, and two SDP relaxations of box-constrained quadratic programming problems, qpG11 and qpG51); {\it(3)} Two infeasible SDPs (infp1 and infd1).

%\begin{enumerate}
%  \item Four small and medium-sized SDPs (two Lov\'asz $\vartheta$ number problems, theta1 and theta2, and two quadratic assignment problems, qap5 and qap9);
%  \item Four large-scale sparse SDPs (two max-cut problems, maxG11 and maxG32, and two SDP relaxations of box-constrained quadratic programming problems, qpG11 and qpG51);
%  \item Two infeasible SDPs (infp1 and infd1).
%\end{enumerate}

Table~\ref{T:SparseStatistic} reports the problem dimensions and some chordal decomposition details. %Figure~\ref{F:SparsityPattern} illustrates the aggregate sparsity patterns of the large-scale sparse SDPs.
The performance of our self-dual method is compared to that of the interior-point solver SeDuMi \citep{sturm1999using}, of the first-order solver SCS~\citep{scs}, and of the primal and dual methods in CDCS~\citep{zheng2016fast}. SparseCoLO~\citep{fujisawa2009user} was used as a preprocessor for SeDuMi. The solution returned by SeDuMi is of high accuracy, so we can assess the quality of the solution computed by CDCS. Instead, SCS is a high-performance first-order solver for general conic programs, so we can assess the advantages of chordal decomposition. Note that SeDuMi should not be compared to the other solvers on CPU time, because the latter only aim to achieve moderate accuracy.
In all cases we set the termination tolerance for CDCS and SCS to $\epsilon_\mathrm{tol} = 10^{-4}$ with a maximum of $2000$ iterations, while we used SeDuMi's default options. All experiments were carried out on a computer with a 2.8 GHz Intel(R) Core(TM) i7 CPU and 8GB of RAM.

Our numerical results are summarized in Tables~\ref{T:ResultsSmall}--\ref{T:ResultsIteration}. In all feasible cases, the objective value returned by our self-dual algorithm was within 0.6\% of the (accurate) value returned by SeDuMi. The small and medium-sized dense SDPs were solved in comparable CPU time by all solvers (Table~\ref{T:ResultsSmall}). For the four large-scale sparse SDPs, our self-dual method was faster than either SeDuMi or SCS (Table~\ref{T:ResultsLargeScale}). As expected, problems with smaller maximum clique size, such as maxG11, maxG32, and qpG11, were solved more efficiently (less than $100\,\rm s$ using our self-dual algorithm).
Note that the conversion techniques in SparseCoLO can give speedups in some cases, but the failure to solve the problem qpG51---due to memory overflow caused by the large number of consensus constraints in the converted problem---highlights the drawbacks.

\begin{table}%[!h]
\centering
\renewcommand\arraystretch{0.68}
\caption{Results for two infeasible SDPs}% in SDPLIB}
\label{T:ResultsInfeasible}
\begin{tabular}{m{0.2cm} r m{0.95cm}<{\centering} m{1.2cm}<{\centering} m{0.95cm}<{\centering} m{0.95cm}<{\centering}}
\toprule[1pt]
& & SeDuMi &\begin{tabular}[x]{@{}c@{}}SparseCoLO\\+SeDuMi\end{tabular}
& SCS &\begin{tabular}[x]{@{}c@{}}Self-dual \end{tabular}\\
\midrule
\multirow{4}{*}{infp1} & Total time (s) & 0.063 & 0.083 & 0.062  & 0.18\\
& Pre- time (s) & 0 & 0.010 & 0.016  &  0.010\\
& Iterations    & 2 & 2 & 20 & 104 \\
& Status        & Infeasible & Infeasible & Infeasible   & Infeasible \\
\midrule
\multirow{4}{*}{infd1}& Total time (s)           & 0.125 & 0.140 & 0.050 & 0.144 \\
& Pre- time (s) & 0 & 0.009 & 0.013  & 0.009\\
& Iterations    & 4 & 4 & 40  & 90 \\
& Status        & Infeasible & Infeasible & Infeasible   & Infeasible \\
\bottomrule[1pt]
\end{tabular}
\end{table}

\begin{table}[t]
\centering
\renewcommand\arraystretch{0.68}
\caption{CPU time per iteration (s)}% for some SDPs in SDPLIB}
\label{T:ResultsIteration}
\begin{tabular}{c c c c c}
\toprule[1pt]
 & SCS &\begin{tabular}[x]{@{}c@{}}CDCS\\(primal)\end{tabular}
&\begin{tabular}[x]{@{}c@{}}CDCS\\(dual)\end{tabular}
&\begin{tabular}[x]{@{}c@{}}Self-dual \end{tabular}\\
\midrule
theta1 & $6\times 10^{-4}$ & $ 2.3 \times 10^{-3}$ &  $2.2 \times 10^{-3}$ &  $2.3\times 10^{-3}$\\
theta2 & $ 1.8 \times 10^{-3}$&$ 5.1 \times 10^{-3}$ &$ 4.7 \times 10^{-3}$ &$ 5.5 \times 10^{-3}$ \\
qap5 & $ 1.2 \times 10^{-3}$ &$ 2.6 \times 10^{-3}$ &$ 2.2 \times 10^{-3}$ &$ 1.9 \times 10^{-3}$ \\
qap9 & $ 1.5 \times 10^{-3}$ & $ 3.6 \times 10^{-3}$ & $ 3.7 \times 10^{-3}$ &$ 4.2 \times 10^{-3}$  \\
\midrule
maxG11 & $ 0.086$ &$ 0.094$ &$ 0.084$ & $ 0.077$\\
maxG32 & 1.243 & 0.260 & 0.231 & 0.209 \\
qpG11 & 0.557 & 0.198 & 0.132& 0.093\\
qpG51 & 1.144 & 0.808 & 0.811& 0.741\\
\bottomrule[1pt]
\end{tabular}
\end{table}

As shown in Table~\ref{T:ResultsInfeasible}, our self-dual algorithm successfully detects infeasible problems, while our previous first-order methods (CDCS-primal and CDCS-dual) do not have this ability. Finally, Table~\ref{T:ResultsIteration} lists the average CPU time per iteration for the first-order algorithms. When comparing the results, it should be kept in mind that our codes are written in MATLAB, while SCS is implemented in C. Nevertheless, we still see that our self-dual algorithm is faster than SCS for the large-scale sparse SDPs (maxG11, maxG32, qpG11 and qpG51), which is expected since the conic projection step is more efficient with smaller PSD cones.

%\balance
\section{Conclusion} \label{se:conclusion}

In this paper, we formulated the homogeneous self-dual embedding of a primal-dual pair of sparse SDPs whose conic constraints are decomposed using chordal decomposition techniques, thereby extending the conversion methods developed in previous work by the authors~\citep{zheng2016fast}. We also showed how the special structure of our homogeneous self-dual formulation can be exploited to develop an efficient ADMM algorithm, which we implemented in the conic solver CDCS.
Our numerical simulations on some benchmark problems from the library SDPLIB show that our self-dual algorithm can give speedups compared to interior-point solvers such as SeDuMi---even when chordal sparsity is exploited using SparseCoLO---and also compared to the state-of-the-art first-order solver SCS. %Future work will exploit chordal sparsity in other areas, such as sparse SDPs from structured controller design~\cite{zheng2016chordal}.%, are also interesting open questions

%Since the current implementation of our algorithms is sequential, but many steps can be carried out in parallel, further computational gains may be achieved by developing our solver CDCS to take full advantage of distributed computing architectures.
%One interesting future work is to apply first-order self-dual embedding algorithms to sparse optimization problems in other areas, such as sparse SDPs from sum-of-squares programming.%, are also interesting open questions.

%\small

\bibliography{ifacconf}             % bib file to produce the bibliography

\begin{thebibliography}{21}
\providecommand{\natexlab}[1]{#1}
\providecommand{\url}[1]{\texttt{#1}}
\providecommand{\urlprefix}{URL }
\expandafter\ifx\csname urlstyle\endcsname\relax
  \providecommand{\doi}[1]{doi:\discretionary{}{}{}#1}\else
  \providecommand{\doi}{doi:\discretionary{}{}{}\begingroup
  \urlstyle{rm}\Url}\fi

\bibitem[{Agler et~al.(1988)Agler, Helton, McCullough, and
  Rodman}]{agler1988positive}
Agler, J., Helton, W., McCullough, S., and Rodman, L. (1988).
\newblock Positive semidefinite matrices with a given sparsity pattern.
\newblock \emph{Linear Alg. Its Appl.}, 107, 101--149.

\bibitem[{Andersen et~al.(2011)Andersen, Dahl, Liu, and
  Vandenberghe}]{andersen2011interior}
Andersen, M., Dahl, J., Liu, Z., and Vandenberghe, L. (2011).
\newblock Interior-point methods for large-scale cone programming.
\newblock \emph{Optim. for machine learning}, 55--83.

\bibitem[{Andersen et~al.(2010)Andersen, Dahl, and
  Vandenberghe}]{andersen2010implementation}
Andersen, M.S., Dahl, J., and Vandenberghe, L. (2010).
\newblock Implementation of nonsymmetric interior-point methods for linear
  optimization over sparse matrix cones.
\newblock \emph{Math. Program. Computation}, 2, 167--201.

\bibitem[{Borchers(1999)}]{borchers1999sdplib}
Borchers, B. (1999).
\newblock {SDPLIB} 1.2, a library of semidefinite programming test problems.
\newblock \emph{Optim. Methods Softw.}, 11(1-4), 683--690.

\bibitem[{Boyd and Vandenberghe(2004)}]{boyd2004convex}
Boyd, S. and Vandenberghe, L. (2004).
\newblock \emph{Convex optimization}.
\newblock Cambridge University Press.

\bibitem[{Fujisawa et~al.(2009)Fujisawa, Kim, Kojima, Okamoto, and
  Yamashita}]{fujisawa2009user}
Fujisawa, K., Kim, S., Kojima, M., Okamoto, Y., and Yamashita, M. (2009).
\newblock User's manual for {SparseCoLO}: Conversion methods for sparse
  conic-form linear optimization problems.
\newblock Technical report, Research Report B-453, Tokyo Institute of
  Technology, Japan.

\bibitem[{Fukuda et~al.(2001)Fukuda, Kojima, Murota, and
  Nakata}]{fukuda2001exploiting}
Fukuda, M., Kojima, M., Murota, K., and Nakata, K. (2001).
\newblock Exploiting sparsity in semidefinite programming via matrix completion
  {I}: General framework.
\newblock \emph{SIAM J. Optim.}, 11(3), 647--674.

\bibitem[{Grone et~al.(1984)Grone, Johnson, S{\'a}, and
  Wolkowicz}]{grone1984positive}
Grone, R., Johnson, C.R., S{\'a}, E.M., and Wolkowicz, H. (1984).
\newblock Positive definite completions of partial hermitian matrices.
\newblock \emph{Linear Alg. Its Appl.}, 58, 109--124.

\bibitem[{Helmberg et~al.(1996)Helmberg, Rendl, Vanderbei, and
  Wolkowicz}]{helmberg1996interior}
Helmberg, C., Rendl, F., Vanderbei, R.J., and Wolkowicz, H. (1996).
\newblock An interior-point method for semidefinite programming.
\newblock \emph{SIAM J. Optim.}, 6(2), 342--361.

\bibitem[{Kalbat and Lavaei(2015)}]{Kalbat2015Fast}
Kalbat, A. and Lavaei, J. (2015).
\newblock A fast distributed algorithm for decomposable semidefinite programs.
\newblock In \emph{IEEE 54th Conf. Decision Control (CDC)}, 1742--1749.

\bibitem[{Kim et~al.(2011)Kim, Kojima, Mevissen, and
  Yamashita}]{kim2011exploiting}
Kim, S., Kojima, M., Mevissen, M., and Yamashita, M. (2011).
\newblock Exploiting sparsity in linear and nonlinear matrix inequalities via
  positive semidefinite matrix completion.
\newblock \emph{Math. Program.}, 129(1), 33--68.

\bibitem[{Madani et~al.(2015)Madani, Kalbat, and Lavaei}]{Madani2015ADMM}
Madani, R., Kalbat, A., and Lavaei, J. (2015).
\newblock {ADMM} for sparse semidefinite programming with applications to
  optimal power flow problem.
\newblock In \emph{IEEE 54th Conf. Decision Control (CDC)}, 5932--5939.

\bibitem[{O'Donoghue et~al.(2016{\natexlab{a}})O'Donoghue, Chu, Parikh, and
  Boyd}]{ODonoghue2016}
O'Donoghue, B., Chu, E., Parikh, N., and Boyd, S. (2016{\natexlab{a}}).
\newblock Conic optimization via operator splitting and homogeneous self-dual
  embedding.
\newblock \emph{J. Optimiz. Theory App.}, 169(3), 1042--1068.

\bibitem[{O'Donoghue et~al.(2016{\natexlab{b}})O'Donoghue, Chu, Parikh, and
  Boyd}]{scs}
O'Donoghue, B., Chu, E., Parikh, N., and Boyd, S. (2016{\natexlab{b}}).
\newblock {SCS}: Splitting conic solver, version 1.2.6.
\newblock \url{https://github.com/cvxgrp/scs}.

\bibitem[{Sturm(1999)}]{sturm1999using}
Sturm, J.F. (1999).
\newblock Using {SeDuMi} 1.02, a {MATLAB} toolbox for optimization over
  symmetric cones.
\newblock \emph{Optim. Methods Softw.}, 11(1-4), 625--653.

\bibitem[{Sun et~al.(2014)Sun, Andersen, and
  Vandenberghe}]{sun2014decomposition}
Sun, Y., Andersen, M.S., and Vandenberghe, L. (2014).
\newblock Decomposition in conic optimization with partially separable
  structure.
\newblock \emph{SIAM J. Optim.}, 24(2), 873--897.

\bibitem[{Vandenberghe and Andersen(2014)}]{vandenberghe2014chordal}
Vandenberghe, L. and Andersen, M.S. (2014).
\newblock Chordal graphs and semidefinite optimization.
\newblock \emph{Found. Trends Optim.}, 1(4), 241--433.

\bibitem[{Ye(2011)}]{ye2011interior}
Ye, Y. (2011).
\newblock \emph{Interior point algorithms: theory and analysis}, volume~44.
\newblock John Wiley \& Sons.

\bibitem[{Ye et~al.(1994)Ye, Todd, and Mizuno}]{ye1994nl}
Ye, Y., Todd, M.J., and Mizuno, S. (1994).
\newblock An $\mathcal{O}\sqrt n l $-iteration homogeneous and self-dual linear
  programming algorithm.
\newblock \emph{Math. Oper. Res.}, 19(1), 53--67.

\bibitem[{Zheng et~al.(2016)Zheng, Fantuzzi, Papachristodoulou, Goulart, and
  Wynn}]{zheng2016fast}
Zheng, Y., Fantuzzi, G., Papachristodoulou, A., Goulart, P., and Wynn, A.
  (2016).
\newblock Fast {ADMM} for semidefinite programs with chordal sparsity.
\newblock \emph{arXiv:1609.06068v2 [math-OC]}.

\bibitem[{Zheng et~al.(2017)Zheng, Fantuzzi, Papachristodoulou, Goulart, and
  Wynn}]{CDCS}
Zheng, Y., Fantuzzi, G., Papachristodoulou, A., Goulart, P., and Wynn, A.
  (2017).
\newblock {CDCS}: Cone decomposition conic solver, version 1.1.
\newblock \url{https://github.com/oxfordcontrol/CDCS}.

\end{thebibliography}

\end{document}